\documentclass[11pt]{amsart}

\usepackage{amsmath,amsfonts,amssymb,amscd}
\begin{document}
\renewcommand{\thefootnote}{\fnsymbol{footnote}}
\pagestyle{plain}

%%%%%%%%%%%%%%%%%%%%%%%%
\title{
{\small \text{Strong K-stability and asymptotic Chow-stability}}
}
\author{Toshiki Mabuchi${}^*$ 
and Yasufumi Nitta${}^{*}{}^{*}$}
\address{Department of Mathematics, Graduate School of Science, 
Osaka University,  Toyonaka, Osaka, 
560-0043 Japan, {\it mabuchi@math.sci.osaka-u.ac.jp}
}
\address{Department of Mathematics, Graduate School of Science, 
Tokyo Institute of Technology,  Megro, Tokyo, 
152-8551 Japan, {\it n$\_$yakkun82@yahoo.co.jp}
}

\maketitle
%$$
%\align
%&x=y
%\tag"(1.1)" \\
%&z=w \tag"(1.2)"
%\endalign
%$$
%%%%%%%%%%%%%%%%%%%%%%%%

\footnotetext{
${}^{*}$Supported 
by JSPS Grant-in-Aid for Scientific Research (A) No. 20244005.\\
 ${}^{*}{}^{*}$Supported 
by JSPS Grant-in-Aid for Young Scientiists (B) No. 23740063.  }

%%%%%%%%%%%%%%%%%%%%%%%%
\abstract
For a  polarized algebraic manifold $(X,L)$, let $T$ be
an algebraic torus in the group $\operatorname{Aut}(X)$ of all holomorphic automorphisms of $X$. 
Then strong relative K-stability (cf.~\cite{M2}) will be shown to imply 
asymptotic relative Chow-stability.  In particular, by taking  $T$ to be trivial, we see that asymptotic Chow-stability
follows from strong K-stability.
\endabstract
%%%%%%%%%%%%%%%%%%%%%%
\section{Introduction}

In this paper, we consider a {\it polarized algebraic manifold $(X,L)$}, i.e., a nonsingular irreducible projective variety $X$, defined over $\Bbb C$, with a very ample line bundle $L$ on $X$.
Let $T$ be an algebraic torus in $\operatorname{Aut}(X)$. Then the main purpose of this paper is to show the following:

\medskip\noindent
{\bf Main Theorem.} 
{\em If $(X,L)$ is strongly K-stable relative to $T$, then $(X,L)$ is asymptotically Chow-stable relative to $T$.}

\section{Relative Chow-stability}

For the maximal compact subgroup $T_c$ of $T$, we put $\frak t_c := \operatorname{Lie}(T_c)$.
For every positive integer $\ell$,  we consider the space $V_{\ell}:= H^0(X,L^{\otimes \ell})$ 
endowed with a Hermitian metric $\rho_{\ell}$ such that the infinitesimal action of $\frak t_c$
on $V_{\ell}$ preserves the metric $\rho_{\ell}$. 
%Take a complex affine line 
%$$
%\Bbb A^1:=\{z \in \Bbb C\},
%$$
Put $\frak t:= \operatorname{Lie}(T)$ and $n:= \dim X$.
Since the infinitesimal action of $\frak t$
 on $X$ lifts to an infinitesimal action of $\frak t$ on $L$, 
we view $\frak t$ as a Lie subalgebra, denoted by $\frak t_{\ell}$, of $\frak{sl}(V_{\ell})$ by taking the traceless part. 
Let $(\frak t_{\ell})_{\Bbb Z}$ be the kernel of the exponential map 
$$
\frak t_{\ell} \owns y \mapsto \exp (2\pi \sqrt{-1}\, y) \in \operatorname{SL}(V_{\ell}).
$$
Let $\frak z_{\ell}$ denote the centralizer of $\frak t_{\ell}$ in $\frak{sl}(V_{\ell})$,
and we consider a symmetric bilinear form $\langle \;, \;\rangle^{}_{\ell}$ on $\frak{sl}(V_{\ell})$ defined by
$$
\; \langle u, v\rangle^{}_{\ell} \; =\; \operatorname{Tr}(u v)/\ell^{n+2},
 \qquad u, v\in \frak{sl}(V_{\ell}),
$$
whose asymptotic limit as $\ell \to \infty$ plays an important role (cf.~\cite{Sz}) 
in the study of relative K-stability for test configurations. 
We now consider the set $\frak t_{\ell}^{\perp}$ 
of all $u \in \frak z_{\ell}$ such that 
$$
\langle u, v\rangle^{}_{\ell} = 0 \qquad  \text{for all $v\in \frak t$.}
$$
By the infinitesimal $\frak t_{\ell}$-action on $V_{\ell}$, we can write the vector space $V_{\ell}$ 
as a direct sum of $\frak t_{\ell}$-eigenspaces
$$
V_{\ell} \; =\; \bigoplus_{k=1}^{m_{\ell}}\; V(\chi_{\ell;k} ),
$$
for mutually distinct additive characters $\chi_{\ell;k} \in 
\operatorname{Hom}((\frak t_{\ell})_{\Bbb Z}, \Bbb Z )$, 
$k = 1,2,\dots, m_{\ell}$, where
 $V(\chi_{\ell;k} )$ denotes the space of all $\sigma\in V_{\ell}$ such that
 $$
 u\,\sigma \; =\; \chi_{\ell;k}(u)\,\sigma
 \qquad \text{for all $u \in (\frak t_{\ell})_{\Bbb Z}$.}
 $$
 Since $T_c$ acts isometrically on $(V_{\ell}, \rho_{\ell})$, the subspaces $V(\chi_{\ell;k} )$ and 
 $V(\chi_{\ell;k'} )$ are orthogonal if $k \neq k'$. For the Lie subalgebra $\frak s_{\ell}$ of 
 $\frak{sl}(V_{\ell})$ defined by
 $$
 \frak s_{\ell}\; =\; \bigoplus_{k=1}^{m_{\ell}}\; \frak{sl} (V(\chi_{\ell;k} )),
 $$
we consider the associated algebraic subgroup $S_{\ell}:= \Pi_{k=1}^{m_{\ell}} 
 \operatorname{SL}(V(\chi_{\ell;k} ))$ of $\operatorname{SL}(V_{\ell})$.
 Let $Z(S_{\ell})$ be the centralizer of $S_{\ell}$ in $\operatorname{SL}(V_{\ell})$.
 Then the Lie algebra $\frak z_{\ell}$ is written as a direct sum of Lie subalgebras
 $$
 \frak z_{\ell} \; =\; \frak z (\frak s_{\ell} ) \oplus \frak s_{\ell}
 $$
where $\frak z (\frak s_{\ell} ):= \operatorname{Lie}(Z(S_{\ell}))$.
For the Lie subalgebra  $\frak t'_{\ell} := \frak t_{\ell}^{\perp}\cap \frak z (\frak s_{\ell} )$ 
of $\frak z (\frak s_{\ell})$, we consider the associated algebraic  subtorus $T'_{\ell}$ of $Z(S_{\ell})$.
Then
$$
T^{\perp}_{\ell} \;:= \; T'_{\ell}\cdot S_{\ell}
$$
is a reductive algebraic subgroup of $\operatorname{SL}(V_{\ell})$ with the Lie algebra 
$\frak t_{\ell}^{\perp}$.
Let $(\frak t_{\ell}^{\perp})_{\Bbb Z}$ denote 
the set of all $u \in (\frak t_{\ell}^{\perp})_{\Bbb Z}$ in the kernel of the 
 exponential map 
$$
\frak z_{\ell} \owns u \mapsto \exp (2\pi \sqrt{-1} u)\in \operatorname{SL}(V_{\ell})
$$ 
such that the circle group $\{ \, \exp (2\pi s\sqrt{-1} u)\,;\, s \in \Bbb R\,\}$ acts isometrically 
on $(V_{\ell}, \rho_{\ell})$.
For each $u \in (\frak t_{\ell}^{\perp})_{\Bbb Z}$, by varying $s \in \Bbb C$, let 
$$
\psi_u : \Bbb G_m \to \operatorname{SL}(V_{\ell}), 
\qquad \exp (2\pi s\sqrt{-1}) \mapsto  \exp (2\pi s\sqrt{-1} u),
$$
be the algebraic one-parameter group generated by $u$, where $\Bbb G_m$ denotes the 
$1$-dimensional algebraic torus $\Bbb C^*$.  
Let $X_{\ell}$ be the image of $X$ 
under the Kodaira embedding
$$
\Phi_{\ell}\; :\; X \to \Bbb P^*(V_{\ell})
$$
associated to the complete linear system $|L^{\otimes \ell }|$ on $X$.
For the degree $d_{\ell}$ of $X_{\ell}$ in $\Bbb P^*(V_{\ell})$, 
we put $W^*_{\ell} := \{\operatorname{Sym}^{d_{\ell}}(V_{\ell}^*)\}^{\otimes n+1}$.
Let $\hat{X}_{\ell} \in  W^*_{\ell}$
be the Chow form for the irreducible reduced algebraic cycle $X_{\ell}$ on $\Bbb P^*(V_{\ell})$,
so that the associated point $[\hat{X}_{\ell}]$ in $\Bbb P^*(W_{\ell})$ is the Chow point for $X_{\ell}$.
Then the action of $T^{\perp}_{\ell}$ on $V_{\ell}$ induces an action of $T^{\perp}_{\ell}$
on $W^*_{\ell}$ and also on $\Bbb P^*(W_{\ell})$.

\medskip\noindent
{\em Definition\/ $2.1$.} (1) $(X,L^{\otimes\ell})$ is called {\it Chow-stable relative to $T$}, if the following conditions are satisfied:

\noindent
(a) \quad The isotropy subgroup of $T^{\perp}_{\ell}$ at 
$[\hat{X}_{\ell}]$ is finite;

\noindent
 (b) \quad The orbit $T^{\perp}_{\ell} \cdot \hat{X}_{\ell}$
in $W_{\ell}^*$ is closed.

\smallskip\noindent
(2)  $(X,L)$ is called {\it asymptotically Chow-stable relative to $T$}, if there exists a positive integer $\ell_0$ 
such that $(X,L^{\otimes \ell})$ are Chow-stable relative to $T$ for all positive integers $\ell$ 
satisfying $\ell \geq \ell_0$.

\section{Test configurations}

Let $u\in (\frak t_{\ell}^{\perp})_{\Bbb Z}$.
For the complex affine line $\Bbb A^1:=\{z \in \Bbb C\}$, we
consider the algebraic subvariety $\mathcal{X}^u$  
of $\Bbb A^1\times \Bbb P^*(V_{\ell})$ obtained as the closure of 
$$
\bigcup_{t\in \Bbb C^*}\; \{t\}\times \psi_{u}(t) X_{\ell}
$$
in $\Bbb A^1\times \Bbb P^*(V_{\ell})$, where $\operatorname{SL}(V_{\ell})$ 
acts naturally on the set $\Bbb P^*(V_{\ell})$ of all hyperplanes in $V_{\ell}$ passing through the origin.
We now put $\mathcal{L}^u:= \operatorname{pr}_2^* \mathcal{O}_{\Bbb P^*(V_{\ell})}(1)$,
where $\operatorname{pr}_2: \mathcal{X}^u \to \Bbb P^*(V_{\ell})$ is the restriction to $\mathcal{X}^u$
of the projection to the second factor: $\Bbb A^1\times \Bbb P^*(V_{\ell})\to\Bbb P^*(V_{\ell})$.
The triple
$$
\mu \; =\; (\mathcal{X}^u, \mathcal{L}^u, \psi_u ), 
$$
is called a {\it test configuration for $(X,L)$ generated by $u$},  where we call $\ell$ the 
{\it exponent} of the test configuration $\mu$.
If $u = 0$, then $\mu$ is called {\it trivial}.

\medskip
For $\mu$ as above, taking the fiber $\mathcal{X}^u_0$ of $\mathcal{X}^u$ over the origin in $\Bbb A^1$, 
we consider the Chow weight $q_{\ell}(u)$ for $\mathcal{X}^u_0$ sitting in $\{0\}\times \Bbb P^*(V_{\ell})
\; (\cong \Bbb P^*(V_{\ell}))$, i.e., the weight at $\hat{\mathcal{X}}^u_0$ of the $\Bbb G_m$-action 
induced by $\psi_u$, where $\hat{\mathcal{X}}^u_0\in W_{\ell}^*$
denotes the Chow form for $\mathcal{X}^u_0$ viewed as an algebraic cycle on $\Bbb P^*(V_{\ell})$.

\medskip\noindent
{\em Definition\/ $3.1$.} (1) $(X,L^{\otimes\ell})$ is called {\it weakly Chow-stable relative to $T$}, if $q_{\ell}(u) < 0$ for all 
$0 \neq u \in (\frak t_{\ell}^{\perp})_{\Bbb Z}$.

\smallskip\noindent
(2) $(X,L)$ is called {\it asymptotically weakly Chow-stable relative to $T$}, if there exists a positive integer $\ell_0$ 
such that $(X,L^{\otimes \ell})$ is weakly Chow-stable relative to $T$ for all positive integers $\ell$ 
satisfying $\ell \geq \ell_0$.

\medskip\noindent
{\em Remark\/ $3.2$.} (1) If  $(X,L^{\otimes\ell})$ is weakly Chow-stable relative to $T$, then by  \cite{M3}, Theorem 3.2, the orbit $T^{\perp}_{\ell}\cdot \hat{X}_{\ell}$ is closed in $W_{\ell}^*$.

\smallskip\noindent
(2) If $(X,L^{\otimes\ell})$ is Chow-stable relative to $T$, then by the Hilbert-Mumford  stability 
criterion, $(X,L^{\otimes\ell})$ is 
weakly Chow-stable relative to $T$.

\section{Strong relative K-stability}

For materials in this section, see  \cite{M} and  \cite{M2}.
For a maximal algebraic torus $\bar{T}$ in $\operatorname{Aut}(X)$
containing $T$,
we fix a Hermitian metric $h$ for $L$ such that $\omega := c_1(L;h)$ is a K\"ahler form preserved by the action 
of the maximal compact subgroup $\bar{T}_c$ of $\bar{T}$.   Then each $V_{\ell}$  admits 
a Hermitian structure $\rho_{\ell}$  preserved by the $T_c$-action
such that
$$
\langle\sigma, \tau\rangle_{\rho_{\ell}}\; :=\; \int_X \,(\sigma, \tau )_{h}\,\omega^n,
\qquad \sigma, \tau\in V_{\ell},
$$
where $(\sigma, \tau )_{h}$ is the pointwise Hermitian pairing of $\sigma$ and $\tau$ on $X$ by
 the Hermitian metric 
$h^{\otimes \ell}$.
In this section, following \cite{M}, we explain how we define the Donaldson-Futaki invariant $F_1$ for a sequence
of test configurations 
$$
\mu_j \; =\; (\mathcal{X}^{u_j}, \mathcal{L}^{u_j}, \psi_{u_j}), \qquad j = 1,2,\dots,
$$
generated by $u_j \in (\frak t_{\ell_j}^{\perp})_{\Bbb Z}$,
where the positive integer $\ell_j$, called the exponent of $\mu_j$, is required to
satisfy $\ell_j \to +\infty$ as $j \to \infty$.
Let $\mathcal{M}$ be the set of all such sequences $\{\mu_j\}$.
For the image ${X}_{\ell_j}$ of $X$ under the Kodaira embedding 
$$
\Phi_{\ell_j}: X \to \Bbb P^*(V_{\ell_j}),
$$
we consider its associated Chow form $\hat{X}_{\ell_j} \in W^*_{\ell_j}:=
\{\operatorname{Sym}_{}^{d_{{\ell}_j}}(V^*_{\ell_j})\}^{\otimes n+1}$. Let $b_{j,\alpha}$, $\alpha = 1,2,\dots, N_{\ell_j}$, be the weights of the $\Bbb G_m$-action
on $V_{\ell_j}$ induced by $\psi_{u_j}$. We then define the norms $\|\mu_j\|^{}_{1}$ and $\|\mu_j\|^{}_{\infty}$ by
$$
\begin{cases}
&\|\mu_j\|^{}_{1} := \; \Sigma_{\alpha=1}^{N_{\ell_j}} \, |b_{j,\alpha}|/\ell_j^{n+1},\\
&\|\mu_j\|^{}_{\infty} := \; \max\{\,  |b_{j,1}|, | b_{j,2}|, \dots,  |b_{j,N_{\ell_j}}|\,\}/\ell_j.
\end{cases}
$$
Let $\delta (\mu_j )$ denote $\|\mu_j \|^{}_{\infty}/\|\mu_j\|^{}_1$ or $1$ according as 
$\|\mu_j\|^{}_{\infty} \neq 0$  or $\|\mu_j\|^{}_{\infty} = 0$.
If $\|\mu_j\|^{}_{\infty} \neq 0$, we write $t\in\Bbb R_+$ as 
$t= \exp (s/\|\mu_j\|_{\infty})$ for some $s \in \Bbb R$, while
we require no relations between $s$ and $t$ if $\|\mu_j\|^{}_{\infty}$ vanishes.
Since $\operatorname{SL}(V_{\ell_j})$ acts naturally on $W^*_{\ell_j}$, 
%by writing $t\in\Bbb R_+$ as 
%$t= \exp (s/\|\mu_j\|_{\infty})$ for $s \in \Bbb R$, 
we define a function $f_{u_j }\,=\, f_{u_j} (s)$ in $s$ on $\Bbb R$ by
$$
f_{u_j}(s)\; :=\; \delta (\mu_j ) \ell_j^{-n}\log \|\psi_{u_j} (t)\cdot \hat{X}_{\ell_j}\|_{\operatorname{CH}(\rho^{}_{\ell_j})}, 
\qquad s \in \Bbb R,
\leqno{(4.1)}
$$
where $W^*_{\ell_j }\owns w \mapsto \|w\|_{\operatorname{CH}(\rho_{\ell_j})} \in \Bbb R_{\geq 0}$ is the Chow norm 
for $W^*_{\ell_j}$ (see \cite{Zh}). Taking the derivative
$\dot{f}_{u_j} (s) := df_{u_j}/ds$, we define $F_1(\{\mu_j \})\in \Bbb R \cup \{-\infty\}$ by
$$
F_1 (\{\mu_j\} ) \;:=\; \lim_{s\to -\infty}\{\varliminf_{j\to \infty} \dot{f}_{u_j} (s)\, \}.
$$
{\em Definition \/$4.2$}\, (cf.~\cite{M2}). (1) $(X,L)$ is called {\it strongly K-semistable relative to $T$}, 
if $F_1 (\{\mu_j\}) \leq 0$ for all $\{\mu_j\}\in \mathcal{M}$.

\smallskip\noindent
(2) Let $(X,L)$ be strongly K-semistable relative to $T$. Then $(X,L)$ is called {\it strongly K-stable relative to $T$}, if for every $\{\mu_j\}\in \mathcal{M}$ satisfying $F_1 (\{\mu_j\}) =0$, there exists a $j_0$ such that 
for all $j \geq j_0$, $\mu_j$ is trivial, i.e., $u_j = 0$. 

\medskip
Note that neither strong K-semistability relative to $T$  nor strong K-stability relative to $T$ depends 
on the choice of $\bar{T}$ and $h$ (see \cite{MN}). 
%For each $u\in (\frak t_{\ell}^{\perp})_{\Bbb Z}$, 
%let  
%$$
%F_1 (\mu_u )\; \in \;\Bbb Q
%$$ 
%denote the ordinary Donaldson-Futaki invariant for the test configuration
%$\mu_u = (\mathcal{X}^u, \mathcal{L}^u, \psi_u)$ generated by $u$.

\section{Proof of Main Theorem}

In this section, we consider a polarized algebraic manifold $(X,L)$ which is strongly K-stable relative to $T$. 
The proof is divided into two parts.

\medskip\noindent
{\em Step\/ $1$.}
We shall first show that $(X,L)$ is asymptotically weakly Chow-stable relative to $T$.
Assume the contrary for contradiction.  Then we can find an increasing sequence of positive integer
$\ell_j$, $j =1,2,\dots$, such that 
$$
\ell_j \to +\infty, \qquad \text{ as $j \to \infty$},
$$
and that $(X, L^{\otimes \ell_j})$ is not weakly Chow-stable relative to $T$ for any $j$.
Then by Definition 3.1, to each $j$, we can assign a element $0  \neq u_j \in (\frak t^{\perp}_{\ell_j})_{\Bbb Z}$ 
such that $q^{}_{\ell_j} (u_j ) \geq 0$. Recall that (see for instance \cite{M0}, Appendix I)
$$
%\begin{align*}
q^{}_{\ell_j} (u_j) \; %&
=\; \|\mu_j\|_1^{}\,\ell_j^n \lim_{s\to -\infty}\dot{f}_{u_j} (s).\; 
%\\
%&=\;
%\; (n+1)! c_1 (L)^n [X]\, \{\, \bar{F}_1(\bar{\mu}_j)\ell_j^n \,+ \,
%\bar{F}_2(\bar{\mu}_j)\ell_j^{n-1} \,+\, \dots\, \}, 
%\end{align*}
%where $\bar{F}_{\alpha} (\bar{\mu}_j )$, $\alpha = 1,2,\dots$,  
%are  the Donaldson-Futaki invariants (see \cite{D}), 
%in an ordinary sense, 
%for test configurations $\bar{\mu}_j := (\mathcal{X}^{u_j},\mathcal{L}^{u_j})$.
$$
Since the function $\dot{f}_{u_j} (s)$ is non-decreasing in $s$ for each $j$, it follows that
$$
0\; \leq \; \|\mu_j \|_1^{-1}\, \ell_j^{-n}q^{}_{\ell_j}(u_j) \; \leq \; \dot{f}_{u_j} (s), \qquad 
-\infty < s < +\infty.
$$
Hence $0 \leq \dot{f}_{u_j}(s)$ for each fixed $s \in \Bbb R$. Taking $\varliminf$  as $j \to \infty$, we have 
$$
0\, \leq \, \varliminf_{j \to \infty}  \dot{f}_{u_j} (s),
\leqno{(5.1)}
$$
for every $s\in \Bbb R$.
By taking limit of (5.1) as $s\to -\infty$, we obtain
$$
0 \; \leq \; \lim_{s\to -\infty}\varliminf_{j \to \infty}  \dot{f}_{u_j}(s)\; = \; F_1 (\{\mu_j \}).
$$
Since $(X,L)$ is strongly K-stable relative to $T$,  this inequality implies that $F_1 (\{\mu_j \})$ vanishes.
Again by strong K-stability of $(X,L)$ relative to $T$, there exists a $j_0$ such that
$\mu_j$ are trivial for all $j$ with $j \geq j_0$ in contradiction to $u_j \neq 0$, as required.

\medskip\noindent
{\em Step\/ $2$.} In view of (1) of Remark 3.2, we see from Step 1 above  that the orbit 
$O_{\ell}:=T^{\perp}_{\ell}\cdot \hat{X}_{\ell}$ is closed in $W_{\ell}^*$. Hence $O_{\ell}$ is an affine 
algebraic subset of $W_{\ell}^*$. Since  $O_{\ell}$ is closed in $W_{\ell}^*$, we here observe that:
$$
O_{\ell} \cap \Bbb C\hat{X}_{\ell} \text{ is a finite set},
\leqno{(5.2)}
$$
where $\Bbb C \hat{X}_{\ell}$ is the one-dimensional vector subspace of 
$W_{\ell}^*$ generated by $\hat{X}_{\ell}$.
Consider the identity components $H_{\ell}$ and $H'_{\ell}$ of the 
isotropy subgroups of the reductive algebraic group
$T^{\perp}_{\ell}$ at the point $\hat{X}_{\ell}$ and $[\hat{X}_{\ell}]$, respectively. 
 Since $\dim H_{\ell} = \dim H'_{\ell}$ by (5.2), it suffices to show that 
 an $\ell_0$ exists 
such  that  $\dim H_{\ell} = 0$ for all $\ell$ with $\ell \geq \ell_0$.
Assume the contrary for contradiction. Then we have an increasing sequence of positive integers $\ell_j$ 
such that 
$$
\dim H_{\ell_j} \; > \; 0,
\qquad j=1,2,\dots,
$$
and that $\ell_j \to +\infty$, as $j \to \infty$.
Since by \cite{Mat} the isotropy subgroup $H_{\ell_j}$ of the reductive algebraic group
$T^{\perp}_{\ell_j}$ at the point $\hat{X}_{\ell_j}$ is a reductive algebraic group,
the group $H_{\ell_j}$ contains a nontrivial algebraic torus $\Bbb G_m$. 
%For the
%maximal connected linear algebraic subgroup $M$ of $\operatorname{Aut}(X)$,
%the infinitesimal action of $\frak m := \operatorname{Lie}(M)$ on $X$
%lifts to an infinitesimal action of $\frak m$ on $L$.
%Hence for each $j$,  we  view $\frak m$ as a Lie subalgebra
% of $\frak{sl}(V_{\ell_j})$.
Since (5.2) allows us to obtain a natural isogeny 
$\iota : H_{\ell_j} \to \bar{H}_{\ell_j}$
 from $H_{\ell_j}$ to an algebraic subgroup 
$\bar{H}_{\ell_j}$ of $\operatorname{Aut}(X)$, 
the image $\bar{H}_{\ell_j}$ also contains a nontrivial algebraic torus $G_j = \Bbb G_m$.
For $\bar{T}$ in Section 4, replacing  $\bar{T}$ by its conjugate in $\operatorname{Aut}(X)$ 
if necessary, 
we may assume that $\bar{T}$ contains $G_j$.
For the maximal compact subgroup $(G_j)_c$ of $G_j$, 
we choose a generator $u_j \neq 0$ for the one-dimensional real Lie subalgebra
$$
\sqrt{-1}\, {(\frak g_j)}_c\; := \; \sqrt{-1}\,\operatorname{Lie}((G_j)_c)
$$
in $ \frak t_{\ell_j}^{\perp}\cap H^0(X, \mathcal{O}(T_X))$ 
such that 
$\exp (2\pi \sqrt{-1} u_j) = \operatorname{id}_X$.
Then for the algebraic group homomorphisms
$\psi^{}_{u_j} : \Bbb G_m \to T^{\perp}_{\ell} \subset \operatorname{SL}(V_{\ell_j})$ generated by $u_j$,
we obtain the associated test configurations 
$$
\mu^{}_{u_j }\; =\; (\mathcal{X}_{}^{u_j}, \mathcal{L}_{}^{u_j}, \psi^{}_{u_j}),
\qquad j=1,2,\dots,
$$
for $(X,L)$ generated by $u_j$.
Let $\beta_j$ be the weight of the $\Bbb G_m$-action by $\psi^{}_{u_j}$ 
at $\hat{X}_{\ell_j}$.
Since $\psi^{}_{u_j}(t)\cdot \hat{X}_{\ell_j} = t^{\beta_j}\hat{X}_{\ell_j} $,
by differentiating the functions $f_{u_j} (s)$ in (4.1) with respect to $s$, we obtain
$$
\dot{f}_{u_j} (s) \; =\; \ell_j^{-n}\beta_j /\|\mu_{u_j}\|^{}_1,
\qquad -\infty < s <+\infty.
\leqno{(5.3)}
$$
Replacing $u_j$ by $v_j := -u_j$, we also have the test configurations 
$$
\mu^{}_{v_j} \; =\; (\mathcal{X}_{}^{v_j}, \mathcal{L}_{}^{v_j}, \psi^{}_{v_j}),
\qquad j = 1,2,\dots,
$$ 
for $(X,L)$ generated by $v_j$. Replace $u_j$ by $v_j$ in (4.1). 
Then by differentiating the functions $f_{v_j} (s)$ with respect to $s$, we obtain
 $$
\dot{f}_{v_j} (s) \; =\; -\,\ell_j^{-n}\beta_j /\|\mu_{v_j}\|^{}_1,
\qquad -\infty < s <+\infty.
\leqno{(5.4)}
$$
Note that $\|\mu_{u_j}\|^{}_1 = \|\mu_{v_j}\|^{}_1$.
The right-hand side of (5.3) 
and the right-hand side of (5.4) are both bounded from above by a positive constant independent of $j$
(see \cite{M}, Section 3). 
Hence, replacing $\{u_j \,;\, j =1,2,\dots \}$ by its subsequence if necessary, we may assume that
$$
\{\,\ell_j^{-n}\beta_j /\|\mu_{u_j}\|^{}_1\;;\; j =1,2,\dots,\,\}
$$
is a convergent sequence. Let $\gamma$ be its limit. Then by (5.3) and (5.4),
$F_1 (\{\mu^{}_{u_j}\}) = \gamma = - F_1 (\{\mu^{}_{v_j}\})$.
Since $(X,L)$ is strongly K-stable relative to $T$, 
the inequalities $F_1 (\{\mu^{}_{u_j}\}) \leq 0$ and $ F_1 (\{\mu^{}_{v_j}\})\leq 0$ hold,
and hence 
$$
\gamma = 0.
$$
Again by strong K-stability of $(X,L)$ relative to $T$, we see that $\mu^{}_{u_j}$ 
are trivial for $j \gg 1$, so that $u_j = 0$ for $j\gg 1$ in contradiction, as required.

%%%%%%%%%%%%%%%%%%%%%%%%%%%%%
%\bigskip\noindent
%{\footnotesize
%{\sc Department of Mathematics}\newline
%{\sc Osaka University} \newline
%{\sc Toyonaka, Osaka, 560-0043}\newline
%{\sc Japan}
%%%%%%%%%%%%%%%%%%%%%%%%%%%%%

%%%%%%%%%%%%%%%%%%%%%%%%%%%%%
%\bigskip\noindent
%{\footnotesize
%{\sc Department of Mathematics}\newline
%{\sc Tokyo Institute of Technology} \newline
%{\sc Meguro, Tokyo, 152-8551}\newline
%{\sc Japan}

\end{document}